\documentclass[11pt]{article}
\usepackage{t1enc}
\usepackage[latin1]{inputenc}
\usepackage[english]{babel}
\usepackage{amsmath,amsthm}
\usepackage{amsfonts}
\usepackage{latexsym}
\usepackage[dvips]{graphicx}
\usepackage{graphicx}
\usepackage[natural]{xcolor}
\usepackage{float}
\usepackage{enumerate}
\usepackage[subnum]{cases}
\usepackage{multirow}
\usepackage{hyperref}
\usepackage{authblk}

\usepackage{listings}
\newtheorem{defi}{Definition}[section]
\newtheorem{thm}[defi]{Theorem}

\newtheorem{lemma}[defi]{Lemma}

\newtheorem{coro}[defi]{Corollary}
\newtheorem{remark}{Remark}
\newtheorem{conj}[defi]{Conjecture}

\newcommand{\TT}{\mathcal{T}}
\newcommand{\TTbar}{\overline{\mathcal{T}}}
\newcommand{\BB}{\mathcal{B}}
\usepackage{lineno}

\date{}
\begin{document}
\baselineskip 0.65cm

\title{\bf Decreasing the mean subtree order by adding $k$ edges}
\vspace{6mm}
\author[1]{Stijn Cambie\thanks{Supported by Internal Funds of KU Leuven (PDM fellowship PDMT1/22/005), UK Research and Innovation Future Leaders Fellowship MR/S016325/1 and the Institute for Basic Science (IBS-R029-C4), stijn.cambie@hotmail.com}}
\affil[1]{Department of Computer Science, KU Leuven Campus Kulak, 8500 Kortrijk, Belgium.}

\author[2]{Guantao Chen\thanks{Partially supported by NSF grant DMS-1855716 and DMS-2154331, gchen@gsu.edu}}
	\affil[2]{Dept.  of Mathematics and Statistics, Georgia State University, Altanta, GA 30303}
	\author[3]{Yanli Hao
		\thanks{The corresponding author, Partially supported by the GSU Provost's Dissertation Fellowship, yhao4@gsu.edu}}
	\affil[3]{School of Mathematics, Georgia Institute of Technology, Atlanta, GA 30332}
	\author[4]{Nizamettin Tokar\thanks{ nizamettintokar@gmail.com}}
	\affil[4]{Department of Mathematics, Usak University, Usak, Turkey 64200}

\date{}

\maketitle

\begin{abstract}
The {\it mean subtree order} of a given graph $G$, denoted $\mu(G)$, is the average number of vertices in a subtree of $G$. 
Let $G$ be a connected graph. 
Chin, Gordon, MacPhee, and Vincent [J. Graph Theory, 89(4): 413-438, 2018] conjectured that if $H$ is a proper spanning supergraph of $G$, then
$\mu(H) > \mu(G)$.  Cameron and Mol [J. Graph Theory, 96(3): 403-413, 2021]  disproved this conjecture by showing that there are infinitely many pairs of graphs
$H$ and $G$ with $H\supset G$, $V(H)=V(G)$ and $|E(H)|= |E(G)|+1$ such that $\mu(H) < \mu(G)$.   They also conjectured that for every positive integer $k$, there exists a pair of graphs $G$ and $H$ with $H\supset G$, $V(H)=V(G)$ and $|E(H)| = |E(G)| +k$ such that $\mu(H) < \mu(G)$. Furthermore, they proposed that $\mu(K_m+nK_1) < \mu(K_{m, n})$ provided $n\gg m$. In this note,  we confirm these two conjectures.

\vskip .2in 

\end{abstract}


\section{Introduction}

Graphs in this paper are simple unless otherwise specified. Let $G$ be a graph with vertex set $V(G)$ and edge set $E(G)$.   The {\it order} of $G$, denoted by $|G|$,  is the number of vertices in $G$, that is, $|G|=|V(G)|$. The {\it complement} of $G$, denoted by $\overline{G}$, is the graph 
on the same vertex set as $G$ such that two distinct vertices of $\overline{G}$ are adjacent if and only if they are not adjacent in $G$. For an edge subset $F\subseteq E(\overline{G})$, denote by $G+F$ the graph obtained from $G$ by adding the edges of $F$.  For a vertex subset $U\subseteq V(G)$, denote by $G-U$ the graph obtained from $G$ by deleting the vertices of $U$ and all edges incident with them.
For any two graphs $G_1,G_2$ with $V(G_1)\cap V(G_2)=\emptyset$, denote by $G_1+G_2$ the graph obtained from $G_1,G_2$ by adding an edge between any two vertices $v_1\in V(G_1)$ and $v_2\in V(G_2)$.

A tree is a graph in which every pair of distinct vertices is connected by exactly one path. A subtree of a graph $G$ is a subgraph of $G$ that is a tree. By convention, the empty graph is not regarded as a subtree of any graph. The {\it mean subtree order} of $G$, denoted $\mu(G)$, 
is the average order of a subtree of $G$.  Jamison~\cite{Jamison1,Jamison2} initiated the study of the mean subtree order in the 1980s, considering only the case that $G$ is a tree. In~\cite{Jamison1}, he proved that  $\mu(T) \ge \frac{n+2}3$ for any tree $T$ of order $n$,   with this minimum achieved if and only if $T$ is a path;  and  $\mu(T)$ could be very close to its order $n$.  Jamison's work on the mean order of subtrees of a tree has received considerable attention~\cite{Haslegrave,Mol,Vince,Wagner,Wang}. At the 2019 Spring Section AMS meeting in Auburn, Jamison presented a survey that provided an overview of the current state of open questions concerning the mean subtree order of a tree, some of which have been resolved~\cite{cambie,luo}.

\begin{figure}[hbt!]
    \centering
    \includegraphics[width=0.2\textwidth]{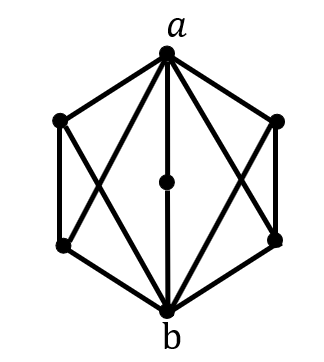}
    \caption{Adding the edge between $a$ and $b$ decreases the mean subtree order}
    \label{ab}
\end{figure}

Recently, Chin, Gordon, MacPhee, and Vincent~\cite{Chin} initiated the study of
subtrees of graphs in general. They believed  that the parameter $\mu$ is monotonic with respect to  
the inclusion relationship of subgraphs. More specifically, they~\cite[Conjecture 7.4]{Chin} conjectured that for any simple connected graph $G$, adding any edge to $G$ will increase the mean subtree order.  Clearly, the truth of this conjecture implies that 
$\mu(K_n)$ is the maximum among all connected simple graphs of order $n$, but it's unknown if  $\mu(K_n)$ is the maximum. Cameron and Mol~\cite{Cameron}  constructed  some counterexamples to this conjecture
 by a computer search.  
 Moreover, they found that the graph depicted in Figure~\ref{ab}  is the  smallest counterexample to this conjecture and   
there are infinitely many graphs $G$ with $xy\in E(\overline{G})$ such that $\mu(G+xy) < \mu(G)$.  In their paper, Cameron and Mol~\cite{Cameron} initially focused on the case of adding a single edge, but they also made the following conjecture regarding adding several edges. 
\begin{conj}\label{ben1}
For every positive integer $k$, there are two connected graphs $G$ and $H$  with $G\subset H$, $V(G)=V(H)$ and $|E(H)\backslash E(G)| =k$ such that  $\mu (H) < \mu (G)$.
\end{conj}

We will confirm Conjecture~\ref{ben1} by proving the following theorem, which will be presented in Section 2. 
 
 \begin{thm}\label{thm-main}For every positive integer $k$, there exist infinitely many pairs of connected graphs $G$ and $H$ with $G\subset H$, $V(G)=V(H)$ and $|E(H)\backslash E(G)| =k$ such that $\mu(H) < \mu(G)$. 
 \end{thm}

In the same paper, Cameron and Mol~\cite{Cameron} also proposed the following conjecture. 
\begin{conj}\label{ben2}
Let $m, n$ be two positive integers.  If $n \gg m$, then we have $\mu(K_m+nK_1) < \mu(K_{m,n})$. 	
\end{conj}

We can derive Conjecture~\ref{ben1} from Conjecture~\ref{ben2}, the proof of which is presented in Section 3, by observing that when $m=2k$, the binomial coefficient ${m \choose 2}$ is divisible by $k$. With $2k-1$ steps, we add $k$ edges in each step, and eventually the mean subtree  order decreases, so it must have decreased in some intermediate step.

\section{Theorem~\ref{thm-main}}

Let $G$ be a graph of order $n$,  and  let $\TT_G$ be the family of subtrees of $G$.  By definition, we have $\mu(G) = (\sum_{T\in \TT_G}|T|)/|\TT_G|$. The {\it density} of $G$ is defined by $\sigma(G) = \mu(G)/n$.  More generally, for any subfamily $\TT\subseteq \TT_G$, we define $\mu(\TT) =(\sum_{T\in \TT}|T|)/|\TT|$ and $\sigma(\TT) = \mu(\TT)/n$. 
 Clearly, $1 \le \mu(G) \le n$ and $0< \sigma(G) \le 1$. 
  
 \subsection{The Construction }
 
 Fix a positive integer $k$. For some integer $m$, let $\{s_n\}_{n\geq m}$ be a sequence of non-negative integers   satisfying: (1) $2s_n\le n-k-1$ for all $n\ge m$; (2) $s_n = o(n)$, i.e., $\lim_{n\to \infty} s_n/n=0$; and (3)
 $2^{s_n} \geq n^2$ for all $n\geq m$. Notice that 
many such sequences exist. Take, for instance, the sequence $\{\lceil 2 \log_2 (n)\rceil\}_{n\ge m}$,  as in \cite{Cameron}, where $m$ is the least positive integer such that $m-2\lceil 2 \log_2 (m)\rceil\ge k+1$.

In the remainder of this paper, we fix $P$ for a path $v_1 v_2 \cdots v_{n-2s_n}$ of order $n-2s_n$. Clearly, $|P|\ge k+1$. Furthermore, let $P^*:= P-\{v_1, \dots, v_{k-1}\}=v_k \cdots v_{n-2s_n}$.

\begin{figure}[hbt!]
    \centering
    \includegraphics[width=0.8\textwidth]{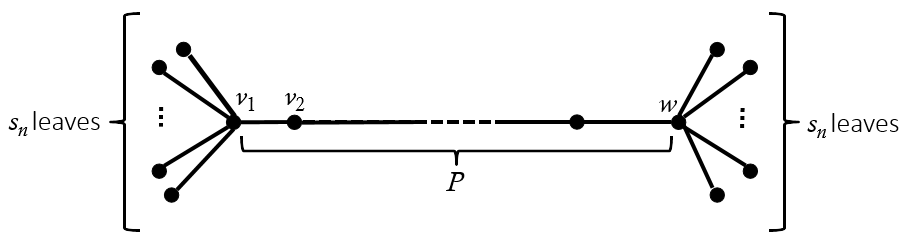}
    \caption{$G_n$}
    \label{G}
\end{figure} 
Let $G_n$ be the graph obtained from the path $P$ by joining $s_n$ leaves to each of the two endpoints $v_1$ and  $w: =v_{n-2s_n}$ of $P$ (see Figure~\ref{G}).  
Let $G_{n,k} :=G_n +\{v_1w, v_2w, \dots, v_kw\}$, that is, $G_{n,k}$ is the graph obtained from $G_n$ by adding 
 $k$ new edges $e_1:= v_1w, e_2:=v_2w, \ldots, e_k:=v_kw$ (see Figure~\ref{G_{n,k}}). 
 \begin{figure}[hbt!]
    \centering
    \includegraphics[width=0.8\textwidth]{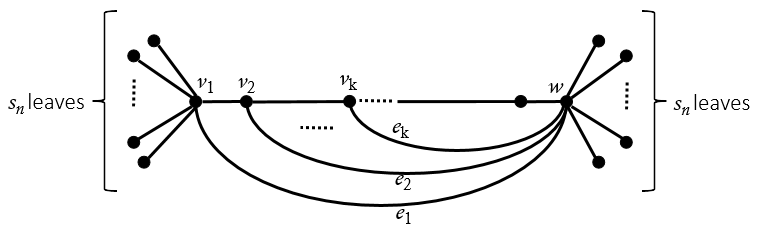}
    \caption{$G_{n,k}$}
    \label{G_{n,k}}
\end{figure}

Let $\TT_{n, k}$ be 
 the family of subtrees of $G_{n,k}$ containing the vertex set $\{v_1, v_k, w\}$ but not containing the path $P^*=v_k \cdots w$.  It is worth noting that $\TT_{n,1}$ is the family of subtrees of $G_{n, 1}$ containing  edge $v_1w$. Note that the graphs $G_n$ and $G_{n, 1}$ defined above are actually the graphs $T_n$ and $G_n$  constructed by Cameron and Mol in \cite{Cameron}, respectively.  From the proof of Theorem 3.1 in~\cite{Cameron}, we obtain the following two results regarding the density of $G_n, G_{n, 1},\TT_{n, 1}$. 
 
 \begin{lemma}\label{thm-Gn} $\lim\limits_{n\to \infty} \sigma(G_n)=1$.
\end{lemma}
  
\begin{lemma} \label{thm-Gn1}
$\lim\limits _{n\rightarrow \infty} \sigma(G_{n, 1})=\lim\limits _{n\rightarrow \infty} \sigma(\TT_{n, 1}) = \frac{2}{3}$.
\end{lemma}

The following two technical results concerning the density of $\TT_{n,k}$ are crucial in the proof of Theorem~\ref{thm-main}. The proofs of these results will be presented in Subsubsection 2.1.1 and Subsubsection 2.1.2, respectively.

 \begin{lemma}\label{thm-Tnk} For any fixed positive integer $k$, $ \lim\limits_{n\rightarrow \infty} \sigma(\TT_{n,k}) = \lim\limits_{n\rightarrow \infty} \sigma(\TT_{n-k+1, 1}).$
 \end{lemma}
 
 \begin{lemma}\label{thm-T2G}
 For any fixed positive integer $k$, $ \lim\limits_{n\rightarrow \infty} \sigma(\TT_{n,k}) = \lim\limits_{n\rightarrow \infty}\sigma(G_{n, k}). $
 \end{lemma}

The combination of Lemma~\ref{thm-Gn1}, Lemma~\ref{thm-Tnk} and Lemma~\ref{thm-T2G} immediately yields the following result.  
\begin{coro}\label{thm-Gnk}
For any fixed positive integer $k$, $\lim\limits _{n\rightarrow \infty} \sigma(G_{n,k}) = \frac{2}{3}$.
\end{coro}
 
 

Combining  Lemma~\ref{thm-Gn} and Corollary~\ref{thm-Gnk}, we have that $\lim\limits_{n\rightarrow \infty}\sigma(G_{n, k})=\frac{2}{3}<1=\lim\limits_{n\rightarrow \infty} \sigma(G_n) $ for any fixed positive integer $k$. By definition,  we gain that $\sigma(G_{n,k}) = \mu(G_{n,k})/|G_{n,k}|$ and $\sigma(G_{n}) = \mu(G_{n})/|G_{n}|$. Since $|G_{n,k}|=|G_{n}|$, it follows that $\mu(G_{n,k}) < \mu(G_n)$ for $n$ sufficiently large, which in turn gives Theorem~\ref{thm-main}.

The following result presented in \cite[page 408, line -2]{Cameron}  will be used in our proof. 

\begin{lemma}\label{tn1}
$|\mathcal{T}_{n, 1}|=2^{2s_n}\cdot\binom{n-2s_n}{2}$.
\end{lemma}


\subsubsection{Proof of Lemma~\ref{thm-Tnk}}

\begin{figure}[hbt!]
    \centering
    \includegraphics[width=0.5\textwidth]{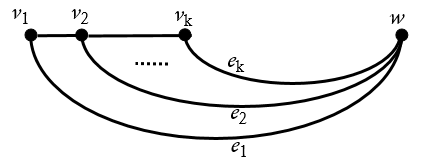}
    \caption{$H$}
    \label{H}
\end{figure}

Let $H$ be the subgraph of $G_{n, k}$ induced by vertex set $\{v_1, \dots, v_k , w\}$ (see Figure~\ref{H}).  Furthermore, set $n_1 = n-k+1$, and let 
 $G_{n_1}^+$ be the graph obtained from $G_{n, k}$ by contracting vertex set $\{v_1, \dots, v_k\}$ into vertex $v_1$ and removing any resulting loops and multiple edges (see Figure~\ref{G_{n_1}^+}).    Clearly, $G_{n_1}^+$ is isomorphic to $G_{n_1,1}$. 
  

\begin{figure}[hbt!]
    \centering
    \includegraphics[width=0.8\textwidth]{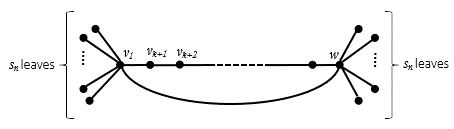}
    \caption{$G_{n_1}^+$}
    \label{G_{n_1}^+}
\end{figure}

Let $T\in\TT_{n, k}$, that is,  $T$ is a subtree of $G_{n,k}$ containing the vertex set  $\{v_1, v_k, w\}$ but not containing the path  
 $P^* =v_k\cdots w$.  Let $T_1$ be the subgraph of $H$ induced by $E(H)\cap E(T)$. Since 
 $T$   does not contain the path $P^*$,  we have that $T_1$ is connected, and so it is a subtree of $H$. 
 Let $T_2$ be the graph obtained from $T$ by contracting vertex set $\{v_1, \dots, v_k\}$ into the  vertex $v_1$ and removing any resulting loops and multiple edges. Since $T_1$ is connected and contains vertex set 
 $\{v_1,v_k,w\}$, it follows that $T_2$ is a subtree of $G_{n_1}^+$ containing edge $v_1w$. So, 
 each $T\in \TT_{n, k}$ corresponds to a unique pair $(T_1, T_2)$ of trees, where $T_1$ is a subtree of $H$ containing vertex set $\{v_1, v_k, w\}$,  
 and $T_2\in \TT_{n_1, 1}$.  We also notice that  $|T|=|T_1|+|T_2|-2$,  where the   $-2$ arises due to the fact that $T_1$ and $T_2$ share exactly 
two vertices $v_1$ and $w$. 
 
  Let $\TT_H'\subseteq \TT_H$ be the family of subtrees of $H$ containing vertex set $\{v_1, v_k, w\}$.  By the corresponding relationship above, 
  we have  $|\TT_{n, k}| = |\TT_H'|\cdot |\TT_{n_1, 1}|$. Hence, we obtain that

\begin{eqnarray*}
    \mu(\TT_{n, k}) &=&\,  \frac{\sum\limits _{T\in \TT_{n,k}}|T|}{|\TT_{n,k}|}= \frac{\sum\limits _{T_1\in \TT_H'} \sum\limits _{T_2\in \TT_{n_1, 1}} \left(|T_1|+|T_2|-2\right)}{|\TT_H'|\cdot |\TT_{n_1, 1}|}\\[1em]
&=&\,\frac{|\TT_H'| \cdot\sum\limits_{T_2\in \TT_{n_1, 1}} |T_2|  + |\TT_{n_1,1}|\cdot\sum\limits _{T_1\in\TT_H'} { |T_1|} - 
    2|\TT_{n_1, 1}|\cdot|\TT_H'|}
    {|\TT_H'|\cdot |\TT_{n_1, 1}|}
    \\[1em]
	&=&\, \mu(\TT_{n_1, 1}) + \mu(\TT_H') -2.
\end{eqnarray*}

Dividing through by $n$, we further gain  that 
$$\sigma(\TT_{n,k})  = \frac{n_1}n\cdot \sigma(T_{n_1, 1})  + 
\frac{k+1}{n}\cdot \sigma(\TT_H')  -\frac{2}{n}.$$
Since $\sigma(\TT_H')$ is always bounded by 1, it follows that $ \lim\limits_{n\rightarrow \infty} \frac{k+1}{n}\cdot \sigma(\TT_H') =0$. Combining this with $\lim\limits_{n\rightarrow \infty} \frac {n_1}n =1$ and $\lim\limits_{n\rightarrow \infty} \frac {2}n =0$, we  get 
$
\lim\limits_{n\rightarrow \infty} \sigma(\TT_{n, k}) = \lim\limits_{n\rightarrow\infty} \sigma(\TT_{n_1, 1})=\frac{2}{3}$ (by Lemma~\ref{thm-Gn1}), which completes the proof of  Lemma~\ref{thm-Tnk}.   \qed

\subsubsection{Proof of Lemma~\ref{thm-T2G}}

Let $\TTbar_{n,k}: =\TT_{G_{n, k}} \backslash \TT_{n, k}$. If $\lim\limits_{n\rightarrow \infty} |\TTbar_{n, k}| /|\TT_{n,k}| =0$, then $\lim\limits_{n\rightarrow \infty}\frac{|\TTbar_{n, k}|}{|\TT_{n,k}| +|\TTbar_{n, k}|}=0$ because $\frac{|\TTbar_{n, k}|}{|\TT_{n,k}| +|\TTbar_{n, k}|}\le |\TTbar_{n, k}| /|\TT_{n,k}|$, and so $\lim\limits_{n\rightarrow \infty}\frac{|\TT_{n, k}|}{|\TT_{n,k}| +|\TTbar_{n, k}|}=1$.  Hence, 
\begin{eqnarray*}
  \lim\limits_{n\rightarrow \infty} \sigma(G_{n, k}) &=&\lim\limits_{n\rightarrow \infty} \frac{\mu(G_{n,k})}{n}=\lim\limits_{n\rightarrow \infty}\frac{1}{n}\cdot \left(\frac{\sum\limits_{T\in \TT_{n,k}}|T|}{|\TT_{n,k}| +|\TTbar_{n, k}|}+\frac{\sum\limits_{T\in \TTbar_{n,k}}|T|}{|\TT_{n,k}| +|\TTbar_{n, k}|} \right)\\&=&\lim\limits_{n\rightarrow \infty}\left(\sigma({\TT_{n, k})}\cdot \frac{|\TT_{n,k}|}{|\TT_{n,k}| +|\TTbar_{n, k}|} + \sigma(\TTbar_{n, k})\cdot\frac{|\TTbar_{n, k}|}{|\TT_{n,k}| +|\TTbar_{n, k}|}\right) =\lim\limits_{n\rightarrow \infty} \sigma({\TT_{n, k}}).
\end{eqnarray*}

Thus, to complete the proof, it suffices to show that 
$\lim\limits_{n\rightarrow \infty} |\TTbar_{n, k}| /|\TT_{n,k}| =0$.  We now define the following two subfamilies of $\TT_{G_{n, k}}$.  
\begin{itemize}
\item $\BB_1=\{ T\in \TT_{G_{n, k}} \ : \ v_1\notin V(T) \mbox{ or } w\notin V(T)\}$; and 
\item $\BB_2 = \{ T\in \TT_{G_{n, k}}\ :  \  T\cap P^* \mbox{ is a path, and $T$ contains $w$}\}$. 
\end{itemize}

Recall that $\TT_{n,k}$ is the family of subtrees of $G_{n, k}$ containing vertex set $\{v_1, v_k, w\}$ and not containing the path $P^*=v_k\cdots w$.  For any $T\in  \TTbar_{n,k}$, by definition, we have the following scenarios: $v_1\notin V(T)$, and so $T\in \BB_1$ in this case; 
$w\notin V(T)$, and so 
$T\in \BB_1$ in this case; 
$v_k\notin V(T)$ and $w\in V(T)$, then $ T\cap P^*$ is a path, and so 
$T\in \BB_2$ in this case; 
 $P^* \subseteq T$, and so 
 $T\in \BB_2$ in this case.   Consequently, $ \TTbar_{n,k}\subseteq  \BB_1\cup  \BB_2 $, which in turn gives that \begin{equation}\label{eq}
 |\TTbar_{n,k}|\le|\BB_1|+|\BB_2|. 
\end{equation}

Let $S_{v_1} $ denote the star centered at $v_1$ with the $s_n$ leaves attached to it and $S_w$ denote the star centered at $w$ with the 
$s_n$ leaves attached to it.  Then $G_{n,k}$ is the union of four subgraphs $S_{v_1}$, $S_w$, $H$, and $P^*$.  
\begin{itemize}
\item Considering the  subtrees of $S_{v_1}$ with at least two vertices and the subtrees of $S_{v_1}$ with a single vertex,  we get  $|\TT_{S_{v_1}}| =(2^{s_n}-1)+(s_n+1)=2^{s_n} +s_n=2^{s_n}+o(2^{s_n})$. 
\item Considering the subtrees of $S_{w}$ with at least two vertices and the subtrees of $S_{w}$ with a single vertex,  we get  $|\TT_{S_w}| =(2^{s_n}-1)+(s_n+1)=2^{s_n} +s_n=2^{s_n}+o(2^{s_n})$. 
\item Considering the subpaths of $P^*$ with at least two vertices and the subpaths of $P^*$ with a single vertex, we get $|\TT_{P^*}| = \binom{|P^*|}{2} + |P^*|= \binom{|P^*|+1}{2}=\binom{n-2s_n -k+2}{2}\le \frac{n^2}{2}$. 
\item The number of subpaths of $P^*$ containing $w$ is bounded above by $|P^*|=n-2s_n-k+1\le n$.
\end{itemize}
Since $s_n=o(n)$, we have the following two inequalities
\begin{eqnarray*}
|\BB_1| & \le & (s_n+ |\TT_H| \cdot |\TT_{P^*}|\cdot|\TT_{S_w}|)+(s_n+|\TT_H| \cdot |\TT_{P^*}|\cdot |\TT_{S_{v_1}} |)\\&\le& 2\left[s_n+|\TT_H|\cdot\left(2^{s_n}+o(2^{s_n})\right)\cdot \frac{n^2}{2}\right] =|\TT_H| \cdot\left(2^{s_n}\cdot n^2+o(2^{s_n}\cdot n^2)\right)\\
|\BB_2| & \le & |\TT_{S_{v_1}}|\cdot |\TT_{S_w}|\cdot |P^*| \cdot |\TT_H| = \left(2^{2s_n}\cdot n + o(2^{2s_n}\cdot n)\right)\cdot |\TT_H|.
\end{eqnarray*}
Recall that $n_1 =n-k+1$.  Applying Lemma~\ref{tn1}, we have  
\begin{eqnarray*}
|\TT_{n, k}|&=& |\TT_H'| \cdot |\TT_{n_1, 1}| = |\TT_H'|\cdot 2^{2s_n} \binom{n_1-2s_n}{2}= |\TT_H'|\cdot2^{2s_n}\cdot \left(\frac{n^2}{2}-o(n^2)\right). 
\end{eqnarray*}

Recall that $2^{s_n} \ge n^2$.  Since $|\TT_H|$ is bounded by a  function of $k$ because $|H| =k+1$, we have the following two inequalities.
\begin{eqnarray*}
\lim_{n\to \infty}\frac{ |\BB_1|}{|\TT_{n, k}|} &=&\lim_{n\to \infty} \frac{|\TT_H|\cdot2^{s_n}\cdot n^2}{|\TT_H'|\cdot2^{2s_n}\cdot \frac{n^2}{2}}=\lim_{n\to \infty} \frac{2|\TT_H|}{|\TT_H'|\cdot 2^{s_n}}=0
\end{eqnarray*}
and
\begin{eqnarray*}
\lim_{n\to \infty}\frac{ |\BB_2|}{|\TT_{n, k}|} &=&\lim_{n\to \infty}\frac{2^{2{s_n}}\cdot n\cdot |\TT_H|}{|\TT_H'|\cdot2^{2s_n}\cdot \frac{n^2}{2}}=\lim_{n\to \infty}\frac{2\cdot |\TT_H|}{|\TT_H'|\cdot n}=0.
\end{eqnarray*}
Hence, we conclude that 
$$\lim\limits_{n \rightarrow \infty}\frac{ |\BB_1|+|\BB_2|}{|\TT_{n, k}|}=0$$ 
Combining this with (\ref{eq}), we have  that $\lim\limits_{n \rightarrow \infty} |\TTbar_{n,k}|/|\TT_{n, k}| =0$, which completes the proof of Lemma~\ref{thm-T2G}.  \qed

\subsection{An Alternative Construction}
The graphs we constructed in order  to prove Theorem~\ref{thm-main}, and the sets of $k$ edges that were added to them, are certainly not the only examples that could be used to prove Theorem~\ref{thm-main}.  For example,  the $k$-edge set $\{v_1w, v_2w, \dots, v_kw\}$ can be replaced by the $k$-edge set $\{v_{1}v_{n-2s_n}, v_{2}v_{n-2s_n-1}$, $\ldots,v_{k}v_{n-2s_n -k+1}\}$. 

Fix a positive integer $k$ and let  $n$ be an integer much larger than $k$.  We follow the notation given in Section 2. Recall that $G_n$ is obtained from a path  $P := v_1 v_2 \cdots v_{n-2s_n}$ by attaching two stars centered at $v_1$ and $v_{n-2s_n}$, and $\lim\limits_{n\to \infty} \sigma(G_n)=1$.   Let $E_k:=\{v_{i_1}v_{j_1}, v_{i_2}v_{j_2}, \ldots, v_{i_k}v_{j_k}\}$  be a set of $k$ edges in $\overline{G_n}$  such that  $1 \le i_1 < j_1 \le i_2 < j_2 \le \dots \le i_k < j_k\le n -2s_n$.  Let $H_{n, k} = G_n + E_k$. 
For convenience, we assume that  $j_\ell -i_\ell$ have the same value, say $p$, for $\ell\in \{1, \dots, k\}$. 

A simple calculation shows that for each path $Q$ of order $q$, we have $\mu(Q) = (q+2)/3$ (See Jamison~\cite{Jamison1}), and so $\lim\limits_{q\to \infty}\sigma(Q) = 1/3$. 
For any non-empty subset $F\subseteq E_k$, we define $\TT_F=\{T\in \TT_{H_{n,k}}:E(T)\cap E_k=F\}$.   For any edge $v_{i_\ell}v_{j_\ell}\in F$, let $e_\ell=v_{i_\ell}v_{j_\ell}$ and $P_\ell =v_{i_\ell}v_{i_\ell+1}\cdots v_{j_\ell}$. Note that every tree $T\in \TT_F$  is 
a union of a subtree of 
$H_{n, k}  -\cup_{e_\ell\in F} (V(P_\ell)\backslash\{v_{i_\ell},v_{j_\ell}\})$ containing $F$ and $\cup_{e_\ell \in F}(E(P_\ell) - E(P_\ell^*))$ for some path  $P_\ell^*\subseteq P_\ell$ containing at least one edge.  Since $|E(P_\ell)| =p$,  the line graph of $P_\ell$  is a path of order $p$. Consequently,  the mean of $|E(P_\ell^*)|$ over subpaths of $P_\ell$ is $(p+2)/3$. Hence,  the mean of  $|E(P_\ell) -E(P_\ell^*)|$ over all subpaths $P_\ell^*$ of $P_\ell$ is $p -(p+2)/3 = 2(p-1)/3$ for each $e_\ell \in F$. Let $s = |F|$. Since every subtree $T\in \TT_F$ has at most $n - s(p-1)$ vertices outside $\cup_{e_\ell\in F} (P_\ell -v_{i_\ell} -v_{j_\ell})$, we get the following inequality. 
\[
\mu(\TT_F) \le  n - s(p-1) + s\cdot \frac{2(p-1)}3  \le  n -\frac{s(p-1)}3.
\] 
By taking $p$ as a linear value of   $n$, say $p =\alpha n$ ($\alpha<\frac{1}{k}$), we get $\sigma(\TT_F) \le 1 - s \alpha/3+s/3n < \sigma(G_n)$ since we assume that $n$ is much larger than $k$. Since   $\TT_{H_{n,k}}=\bigcup_{F\subseteq E_k}\TT_F$, we have $\sigma(H_{n, k}) < \sigma(G_n)$, and so $\mu(H_{n, k}) < \mu(G_n)$. 

\begin{remark}
    The above construction gives an example where we can delete $k$ edges in order in such a way that the mean subtree order increases in every step.
\end{remark}

\section{Proof of Conjecture~\ref{ben2}}

To simplify notation, we let $G:=K_m+nK_1$, where $V(G) = V(K_{m, n})$. Denote by $A$ and $B$ the two color classes of $K_{m, n}$ with $|A|=m$ and $|B|=n$, respectively. For each tree $T\subseteq G$, we have $E(T)\cap E(K_m)=\emptyset$ or $E(T)\cap E(K_m)\ne\emptyset$. This implies that the family of subtrees of $G$ consists of the subtrees of $K_{m,n}$ and the subtrees sharing at least one edge with $K_m$.
For each tree $T\subseteq G$, let $A(T) = V(T)\cap A$ and $B(T) = V(T)\cap B$. Then, $|T|=|A(T)|+|B(T)|$. Furthermore, let $B_2(T)$ and $B_{\geq 2}(T)$ be the sets of vertices $v\in B(T)$ such that $d_T(v) = 2$ and $d_T(v) \geq 2$, respectively. Clearly, $B_2(T)\subseteq B_{\geq 2}(T)\subseteq B(T)$.
We define a subtree $T\in\TT_G$ to be a {\it b-stem} if $B_{\geq 2}(T) = B(T)$, which means that $d_T(v)\ge 2$ for any $v\in B(T)$.  

Let $T$ be a b-stem and assume that $T$ contains $f$ edges in $K_m$. Counting the number of edges in $T$, we obtain
$|E(T)| = f+\sum_{v\in B(T)}d_T(v)$. Since $T$ is a tree, we have $|E(T)| = |T| -1 =  |A(T)| + |B(T)| -1$. 
Therefore, we gain 
\begin{equation}\label{eqn-2}
|B(T)| = |A(T)| -1-\left(f +\sum_{v\in B(T)}(d_T(v) -2)\right).
\end{equation}
Since $T$ is a b-stem, we have $\sum_{v\in B(T)}(d_T(v) -2) \ge 0$, which implies that $|B(T)|\leq |A(T)| -1\leq m-1$. Thus, $|T|=2 |A(T)|-\left(1+f+\sum_{v\in B(T)}(d_T(v) -2)\right) \le 2|A(T)| -1$. It follows that  a b-stem $T\in \TT_G$ is the {\it max b-stem}, i.e., the b-stem with the maximum order among all b-stems in $\TT_G$, if and only if $A(T) = A$, $E(T)\cap E(K_m)=\emptyset$, and $B_2(T)=B_{\geq 2}(T)$. This is equivalent to saying that  $T$ is a max b-stem if and only if $|A(T)|=m$ and $|B(T)|=m-1$.

The {b-stem} of a tree $T\subset G$ is the subgraph induced by $A(T)\cup B_{\geq 2}(T)$, and it is a subtree in $\TT_G$.  It is worth noting that the b-stem of every subtree $T\subset G$ exists, except for the case when $T$ is a tree with only one vertex belonging to $B$.  Conversely, given a b-stem $T_0$, a tree $T\subset G$ contains $T_0$ as its b-stem if and only if $T_0\subseteq T$, $A(T)=A(T_0)$, and $B(T)\backslash B(T_0)$ is a set of vertices with degree 1 in $T$. Equivalently, $T$ can be obtained from $T_0$ by adding vertices in $B(T)\backslash B(T_0)$ as leaves. So, there are exactly $(|A(T_0)|+1)^{n-|B(T_0)|}$ trees containing $T_0$ as their b-stem.

For two non-negative integers $a, b$, where $a \ge b+1 \ge 1$,
let $\TT_G(a,b)$ (resp. $\TT_{K_{m,n}}(a,b)$) be the family of subtrees in $\TT_G$ (resp. $\TT_{K_{m,n}}$) whose b-stems $T_0$ satisfy $|A(T_0)|=a$ and $|B(T_0)|=b$. For any $A_0\subseteq A$ and $B_0\subseteq B$, let $f_G(A_0,B_0)$ (resp. $f_{K_{m,n}}(A_0,B_0)$) denote the number of b-stems $T_0$ spanned by $A_0\cup B_0$; that is, $A(T_0)=A_0$ and $B_{\ge 2}(T_0)=B_0$. 
Clearly, $f_G(A_0,B_0)$ and $f_{K_{m,n}}(A_0,B_0)$ depend only on $|A_0|$ and $|B_0|$, so we can denote them by $f_G(|A_0|,|B_0|)$ and $f_{K_{m,n}}(|A_0|,|B_0|)$, respectively. 
By counting, we have $|\TT_G(a,b)|=\binom{m}{a}\cdot\binom{n}{b}\cdot f_G(a,b)\cdot (a+1)^{n-b}$ and $|\TT_{K_{m,n}}(a,b)|=\binom{m}{a}\cdot\binom{n}{b}\cdot f_{K_{m,n}}(a,b)\cdot (a+1)^{n-b}$, due to the fact that there are $\binom{m}{a}$ ways to pick an $a$-set in $A$ and $\binom{n}{b}$ ways to pick a $b$-set in $B$. Since $a\leq m$ and $b\leq m-1$, there exist positive numbers $c_1$ and $c_2$ that depend only on $m$, such that
\begin{equation}\label{eqn-3}
	c_1n^{b} (a+1)^{n-b}\le |\TT_G(a,b)| \le c_2n^{b}(a+1)^{n-b}
\end{equation}
 Note that if $(a,b)\ne (m,m-1)$, then we have $b\le m-2$. Applying inequality (\ref{eqn-3}), we get $|\cup_{(a,b)\ne(m,m-1)}\TT_G(a,b)| \le c_3|\TT_G(m,m-1)|/n$ for some constant $c_3 > 0$ depending only on $m$. 
 
Given a  b-stem $T_0$ with $|A(T_0)|=a$ and $|B(T_0)|=b$, let $T$ be a tree chosen uniformly at random from $\TT_G$ (resp. $\TT_{K_{m,n}}$) that contains $T_0$ as its b-stem. Then, the probability of a vertex $v\in B\backslash B(T_0)$ in $T$ is $\frac{a}{a+1}$. This shows that the mean order of trees containing $T_0$ as their b-stem is $(n-b)\frac{a}{a+1}+a+b$, denoted by $\mu(a,b)$. Note that $\sum_{T\in \TT_G(a,b)}|T|=\mu(a,b)\cdot |\TT_G(a,b)|$ and $\sum_{T\in \TT_{K_{m,n}}(a,b)}|T|=\mu(a,b)\cdot |\TT_{K_{m,n}}(a,b)|$. Assume that $T_0$ has $f$ edges in $K_m$, and set $c=\sum_{v\in B(T_0)}(d_{T_0}(v)-2)$. Using (\ref{eqn-2}), we have $b=a-(1+f+c)$. Hence, $\mu(a,b)=\frac{(n+2+a)\cdot a}{a+1}-\frac{1+f+c}{a+1}$, which reaches its maximum value when $a=m$ and $f=c=0$, i.e., when $T_0$ is a max b-stem.
We then have:

\begin{eqnarray*}
\mu(G)&=&\frac{\mu(m,m-1)|\TT_G(m,m-1)|+\sum_{(a,b)\ne (m,m-1)}\mu(a,b)|\TT_G(a,b)|+n}{|\TT_G(m,m-1)|+\sum_{(a,b)\ne (m,m-1)}|\TT_G(a,b)|+n},
\end{eqnarray*}

\begin{eqnarray*}
\mu(K_{m,n})&=&\frac{\mu(m,m-1)|\TT_{K_{m,n}}(m,m-1)|+\sum_{(a,b)\ne (m,m-1)}\mu(a,b)|\TT_{K_{m,n}}(a,b)|+n}{|\TT_{K_{m,n}}(m,m-1)|+\sum_{(a,b)\ne (m,m-1)}|\TT_{K_{m,n}}(a,b)|+n},
\end{eqnarray*}
where $n$ denotes the number of subtrees with a single vertex in $B$.

Note that $|\TT_G(a,b)|\ge |\TT_{K_{m,n}}(a,b)|$, with equality holding if and only if $a=b-1$, and so in particular when $(a,b)=(m,m-1).$
We have derived before that $0 <\mu(a,b)<\mu(m,m-1)$ when $(a,b)\ne (m,m-1).$
Using the inequality $|\cup_{(a,b)\ne(m,m-1)}\TT_G(a,b)|\le c_3|\TT_G(m,m-1)|/n$, we conclude that $\mu(G)> \frac{n}{n+c_3} \mu(m,m-1)> \max_{ (a,b)\ne (m,m-1)} \mu(a,b)$ for $n$ sufficiently large (for fixed $m$).

Since $\mu(K_{m, n})$ is the average of the same terms, as well as some additional terms of the form $\mu(a,b)$, which are smaller than $\mu(G)$, we conclude that $\mu(G) < \mu(K_{m, n})$.
This completes the proof. \qed

 \section*{Acknowledgments}
We would like to express our sincere gratitude to the  anonymous referees for their valuable comments and suggestions that improved this manuscript.

\bibliographystyle{plain}

\bibliography{LA-Ref.bib}

\end{document}